\documentclass[12pt,notitlepage]{article}
\usepackage{epsfig, xypic,eucal,amsmath, amssymb, amsthm,latexsym}

\newtheorem{nle}{Lemma}

\newtheorem{te}{Theorem}

\newtheorem{cl}{Corollary}

\newtheorem{defi}{Definition}

\newcommand{\frb}{\mathfrak{b}}

\newcommand{\frg}{\mathfrak{g}}                                                 
\newcommand{\frh}{\mathfrak{h}}

\newcommand{\frk}{\mathfrak{k}}

\newcommand{\frp}{\mathfrak{p}}                                                 
\newcommand{\frqq}{\mathfrak{q}}

\begin{document}

\title{Transformation formulas in Quantum Cohomology}
\author{Prakash Belkale\\
Department of Mathematics\\
University of Utah\\
155 S 1400 E JWB 233\\
Salt Lake City, Utah 84112-0090\\
      belkale@math.utah.edu                         }
\maketitle

\begin{abstract} This article discusses equality of certain Gromov-Witten numbers of  $G/B$'s and $G/P$'s which are suggested by some problems in representation theory. We will see that this translates into an 'action of the center' which has many applications.  
\end{abstract}

\section{Introduction}

It is known \cite{wood},\cite{be} that the problem of determining the conditions on conjugacy classes $\bar{A_1},\dots,\bar{A_s}$ in $SU(n)$ so that
these lift to elements $A_1,\dots,A_s\in SU(n)$ with $A_1A_2,\dots A_s = 1$
is controlled by Quantum Schubert calculus of Grassmannians. Teleman and
Woodward have in a recent preprint generalized this to an arbitrary
simple, simply connected compact group $K$. If $G$ is the complex simple subgroup (whose real points are $K$), then the role played by the Grassmanians is
replaced by the homogeneous spaces $G/P$ for P a maximal parabolic subgroup.

In the case of $SU(n)$ (and similarly for $K$) there is a natural 'action'
of center of $K$ on the representation theory side, namely if $c_1,\dots,c_s$
are central elements with $c_1c_2\dots c_s=1$, then these act on the set of
conjugacy classes $\bar{A_1},\dots,\bar{A_s}$ in $SU(n)$  lift-able to elements $A_1,\dots,A_s\in SU(n)$ with $A_1A_2\dots A_s = 1$, the action being just
multiplying $\bar{A_i}$ by $c_i$. This action is well defined on the level
of conjugacy classes because the $c_i$ are central.

This suggests a natural transformation property of Gromov-Witten 
numbers of the Grassmannians under the action of the center. This property
was proved in \cite{wood} as a consequence of the known description of
Quantum Schubert calculus \cite{bert}.

The aim of this article is twofold. The first aim is to prove the transformation formulas geometrically and in complete generality (for any simple simply connected complex Lie group). The second is to show that these formulas determine
quantum Schubert calculus in the case of Grassmannians (Bertram's Schubert Calculus). We also give a strengthening in the case of Grassmannians of a theorem of
Fulton and Woodward on the lowest power of $q$ appearing in a (quantum)product of Schubert classes in $G/P$ where $P$ is a {\em maximal} parabolic subgroup.

Let us now describe these transformation formulas. See section 1 for the notation. Let $G$ be a simple, semi-simple  complex algebraic group. We first construct a map
$\phi: C\to W$ where $C$ is the center of $G$ and $W$ the Weyl group.
Let $Z$ be a homology class of $G/P$ where $P$ is an arbitrary parabolic
subgroup (not necessarily maximal). Let $c_1,\dots,c_s$ be central
elements with product =1. Let $w_1,\dots,w_s$ be elements of (a suitable
right quotient of $W$), then the transformation formulas take the shape

$$<X_{\phi(c_1)w_1},\dots,X_{\phi(c_s)w_s}>_{Z'}=<X_{w_1},\dots,X_{w_s}>_Z$$

where $Z'$ is a homology class determined by $Z$ and the rest of the data.

It happens that in some cases, $Z'$ is a simpler homology class than $Z$,
for instance $Z'$ could be zero, when $Z$ is not. This allows
for reducing the quantum terms to the classical ones. This program works in the Grassmannian case (and Flag manifold case, but we have not included the details here).

There exist simple simply connected groups with trivial center. In this case
the transformation formulas do not give any information. Maybe there is
an extension of these transformation formulas to the non miniscule case.
Such an extension is not apparent from the representation theory side.

The argument for Grassmannians may go through in the miniscule case. We will
return to this question in a later paper.

One final comment is that even in the classical case of cohomology the
transformation formulas give vanishing statements. For example if
Z' turns out to be negative and $Z=0$, then we get a vanishing statement
of certain intersection numbers.

Many of the results in this paper are new proofs of older results using methods which seem both natural and elementary (to the author). It is perhaps worth pointing out what is essentially new in this paper - the transformation formulas in the (usual) partial flag manifold
case, the exact determination of the lowest order terms in the quantum products of two Schubert cycles in (usual) Grassmannians, the natural extension of the transformation formulas to all groups. It should be noted that the quantum schubert calculus of $G/B$ or of $G/P$ is not known except (when $G=\text{SL}(n)$), and the relations coming from the center might be the first general explicit formula.

A. Buch \cite{buch} has recently given new proofs of Bertram's Quantum Schubert Calculus, using very different methods. 

We thank Aaron Bertram, Patrick Brosnan, Anders Buch, William Fulton, Misha Kapovich, Madhav Nori and Christopher Woodward for several helpful conversations and encouragement.

\section{Some representation theory}

\subsection{Notation}
 We review some basic representation theory in this section. For proofs refer
to Bourbaki\cite{bourbaki}.

Let $G$ be a simple simply connected complex algebraic group. Let $\frg$ be its lie algebra. Let $B$ be a Borel subgroup, $T\subset B$ a maximal torus and
let

$$ \frg = \frh \bigoplus_{\alpha}\frg_{\alpha}.$$

where the $\alpha$'s belong to the subset of roots $R$ in $\frh^{*}$. The set
$R$ is partitioned into the set of positive roots $R^{+}$ and negative roots
$R^{-}$, and the lie algebra of $B$ is 

$$\frb = \frh \bigoplus_{\alpha}\frg_{\alpha}.$$ 

with the $\alpha's$ in $ R^{+}$. Also define $\Delta$ to be the set of simple roots.

The Weyl group W is defined to be $N(T)/T$, where $N(T)$ is the normalizer
of T which acts on $\frh$ and $\frh^*$. If $\alpha$ is a root, we have elements
$w_{\alpha} \in W$, $H_{\alpha} \in \frh$, so that $w_{\alpha}$ acts on
$\frh^*$ by

$$w_{\alpha}(\beta) = \beta -\beta(H_{\alpha})\alpha.$$

and this map preserves the roots, is a reflection and takes $\alpha$ to $-\alpha$.

The action of $w_{\alpha}$ on $\frg$ is given by

$$w_{\alpha}(H) = H -\alpha(H)H_{\alpha}.$$

It is also known that $w_{\alpha}$'s generate $W$.

The affine Weyl group $W_{aff}$ is defined to be the set of automorphisms
of $\frh$ generated by $W$ and translations by $H_{\alpha}$ for $\alpha \in R$.

\subsection{Conjugacy classes}

Let $K$ be a maximal (connected) compact subgroup of $G$. If $\frk$ is the lie
algebra of $K$, then $\frk \bigotimes_{\Bbb{R}} \Bbb{C} = \frg$. Let $T_K = T \cap K$
be the maximal torus in $K$, with lie algebra $i\frh$. The following are standard facts:

\begin{enumerate}

\item $T_K \to K$ induces a surjection on conjugacy classes.

\item $x_1,x_2 \in T_K$ are conjugate in $K$ if and only if  there exists $w \in W$ with $Ad(w)x_1 = x_2$. 

\item let $\text{Exp}: \frh_{K} \to T_K$ be the exponential map from the lie algebra
of $T_K$ to $T_K$. The kernel of this map is $\Gamma(T) =  \Bbb{Z}-$span $\{ 2\pi i 
H_{\alpha}: \alpha \in R\}$. This follows from simply connectedness of $G$.

\item If $t_1, t_2 \in \frh_K$, then \text{Exp}($t_1$) and Exp($t_2$) are conjugate
in $K$ if and only if there exists $w \in W$ with 

$$w(t_1)-t_2 \in \Gamma(T).$$
\end{enumerate}

Putting this all together, if we let $\frh_{\Bbb{R}}$ be the $\Bbb{R}$- span 
of $\{H_{\alpha}: \alpha \in R\}$, then the map
$\frh_{\Bbb{R}} \to T_K$ by $t \mapsto $\text{Exp}($2\pi it$) induces an isomorphism

$\frh_{\Bbb{R}}/ W_{aff} \to $ conjugacy classes in $K$.

\subsection{Fundamental Chamber and Center}

Let $L_{\alpha,k}=\{ x\in \frh_{\Bbb{R}}: \alpha(x) = k\}$. The affine Weyl group
is then the group generated by reflections in $L_{\alpha,k}$ for $k \in \Bbb{Z}$. Finally let $\breve{\alpha}$ be the highest weight for the adjoint representation.

\begin{te}(Fundamental chamber for affine Weyl group)
Let $C= \{x \in \frh_{\Bbb{R}}: \alpha(x)>0$ for $ \alpha \in R^{+}, \breve{\alpha}(x)<1\}$.
\begin{enumerate}
\item $C$ is a connected component of $\frh_{\Bbb{R}} - \cup_{\alpha \in R, k \in \Bbb{Z}} L_{\alpha,k}.$ 
\item If $C'$ is any other component, there is a unique $w \in W_{aff}$ with
$w(C) = C'.$
\item Let $\bar{C}$ be the closure of $C$, then the composite 
$\bar{C}\to \frh_{\Bbb{R}}/ W_{aff} \to $ conjugacy classes in $K$ is, a
homeomorphism.
\end{enumerate}
\end{te}

We now give the description of the center. For this, Let $S =\{x \in \bar{C}:
\alpha(x) \in \Bbb{Z}$, for all $\alpha \in R \}$. Finally, write
$\breve{\alpha} = \sum_{\alpha \in \Delta} n_{\alpha} \alpha$. where $\Delta$ is the set of simple positive roots.

\begin{te}
\begin{enumerate}
\item The map $\bar{C}\to $ conjugacy classes in $G$, takes $S$ to center($K$).
\item Define $x_{\alpha}$ for $\alpha \in \Delta$ by the formula
$\beta(x_{\alpha})= \delta_{\alpha,\beta}$ for $\alpha, \beta \in \Delta$.
Then, $S = \{0\} \cup \{x_{\alpha}: \alpha \in \Delta, n_{\alpha}=1\}$ 
\end{enumerate}
\end{te}

\subsection{ A map from Center to the Weyl group}

If $c$ is in center($G$) (=center($K$)), which corresponds to $x_{\alpha} \in \Delta
$ then we can consider the set $C'=C - x_{\alpha}$. It is easy to see that
this is a connected component of $\frh_{\Bbb{R}} - \cup_{\alpha \in R, k \in \Bbb{Z}} L_{\alpha,k}.$  Therefore we have 

$$C - x_{\alpha} = w_{c}^{-1}(C) + t_{\alpha}$$ 

where $w_c \in W$ and $t_{\alpha} \in \Bbb{Z}$ - span of $H_{\beta}, \beta \in R.$ \\

Define $y(x) = w_{c}(x-x_{\alpha} - t_{\alpha})$. Then we have
the equation

$$x-x_{\alpha} = w_{c}^{-1}(y(x)) + t_{\alpha}$$

and also that $x\in C \implies y(x) \in C$. Hence, $x\in \bar{C} \implies y(x) \in \bar{C}.$ Now put $x= x_{\alpha}$.  Then we get,

$ w_{c}^{-1}(y(x) = -t_{\alpha}$, this gives $y(x)$ is zero in $\frh_{\Bbb{R}}/ W_{aff}$ and by theorem 1 in section 2.3, we get $y(x)=0$ so $t_{\alpha}=0.$

We therefore have:
\begin{te}(Map from center to Weyl group)
Let $c$ be in center($K$), let $\alpha \in \Delta$ be so that 
$exp(2\pi i x_{\alpha}) = c $. Then there exists a $w_c\in W$ so that
the equation $C - x_{\alpha} = w_{c}^{-1}(C)$ holds. Furthermore
the map Center($G$)$ \to W$ is a injective homomorphism of groups.
\end{te}
\begin{proof}
The only part not proved is that  $c\mapsto w_c$ is a homomorphism
of groups. For this, let $c_1, c_2$ be central elements. Let $x_{\alpha},x_{\beta},x_{\gamma}$ correspond to $c_1, c_2$ and $c_1 c_2$ respectively.
It is clear that $x_{\alpha}+x_{\beta}=x_{\gamma} +t$ with $t\in \Bbb{Z}$ - span of $H_{\delta}, \delta\in R.$ 
hence $$C-x_{\gamma}= C - x_{\alpha}-x_{\beta} + t$$

$$ =w_{c_1}^{-1}(C -x_{\beta})+ (w_{c_1}^{-1}(x_{\beta})-x_{\beta}) +t$$

$$= w_{c_1}^{-1}w_{c_2}^{-1}(C) +t_1 +t$$

where $t_1$ and $t$ are in $\Bbb{Z}$ - span of $H_{\delta}, \delta\in R.$
The proof is therefore complete.
\end{proof}

We can describe the element $w_{c}$ more concretely, for this first note
\begin{enumerate}
\item If $x \in C$ then $y(x) \in C$.

\item For $\beta \in R$, $\beta \in R^{+}$ if and only if  $\beta(x) >0$ for any $x \in C.$

\item $w_c(\beta)(y(x)) =\beta(w_c^{-1}(y(x)))= \beta(x) - \beta(x_{\alpha}).$ 

\end{enumerate}

we therefore have the following description of $w_c$.
\begin{nle} In the situation above 
\begin{enumerate}
\item If $\beta(x_{\alpha}) =0$ then $\beta$ and $w_{c}(\beta)$ have the
same 'sign' (simultaneously in $R^{+}$, or in $R^{-}$).
\item If $\beta(x_{\alpha}) =1$ then $\beta$ and $w_{c}(\beta)$ have opposite
'signs'.
\end{enumerate}
 
\end{nle}

\subsection{Parabolics associated to Central Elements}

Fix $c_1$, $c_2$ belonging to center($G$) with $c_1 c_2 =1.$ Let $x_{\alpha},
x_{\beta}$ be their representatives in $\bar{C}.$ 

Define parabolic subgroups $P_1,P_2 \supset B$ by defining their Lie algebras

$$\frp_1 = \frh \bigoplus_{\gamma: \gamma(x_{\alpha})\geq 0}\frg_{\gamma}.$$

$$\frp_2 = \frh \bigoplus_{\gamma: \gamma(x_{\beta})\geq 0}\frg_{\gamma}.$$

And their Levi subgroups $Q_1,Q_2$ by defining their lie algebras

$$\frqq_1 = \frh \bigoplus_{\gamma: \gamma(x_{\alpha})= 0}\frg_{\gamma}.$$

$$\frqq_2 = \frh \bigoplus_{\gamma: \gamma(x_{\beta}) = 0}\frg_{\gamma}.$$

That these are closed subgroups, follows from the $P$'s being standard parabolics and $Q$'s being  the  centralisers of the $c$'s.

These are related due to the relation $c_1 c_2 =1$:

\begin{nle}
$Q_2 = Ad(w_{c_1})(Q_1)$. Or that $\gamma(x_{\alpha}) = 0$ if and only if
$(w_{c_1}(\gamma))(x_{\beta}) = 0$. 
\end{nle}
\begin{proof} Follows from $-x_{\alpha}= w_{c_1}^{-1}(x_{\beta}).$
\end{proof}

\begin{nle}
$(Ad(w_{c_1})(P_1))\cap P_2 =Q_2$ and this is a transverse intersection.
\end{nle}
\begin{proof}
The first statement follows from
$$(w_{c_1}(\gamma))(x_{\beta}) = -\gamma(x_{\alpha}).$$ The transversality statement follows from the same equation (counting dimensions).
\end{proof}

\begin{cl} If $g_1,g_2$ are general elements of $G$, the set
$(w_{c_1} P_1 g_1)\cap (P_2 g_2)$ is non empty.
\end{cl}
\begin{proof} If $g_1 = w_{c_1}^{-1}, g_2 = 1$, then this follows from the lemma above. Then apply standard intersection theory(local).
\end{proof}

\section{Algebraic Geometry preliminaries on $G/B$}

\subsection{Line bundles on $G/P$ and $G/B$}
Let $P$ be a parabolic containing $B$. We then have a natural surjection
$G/B\to G/P$. It is known this induces injections on the Picard Groups.
Our goal here is to recall the standard facts on describing all the line bundles on $G/B$ and on which of these descend to $G/P$.

Let $WL$ = weight lattice of $\frg$. This is the subset of $\frh^*$ spanned by
elements $\omega$ so that $\omega(H_{\alpha})\in \Bbb{Z}$ for all $\alpha \in R$
. It has a $\Bbb{Z}$ basis $\{\omega_{\alpha}:\alpha \in \Delta\}$ where 

$$\omega_{\alpha}(H_{\beta}) =\delta_{\alpha,\beta}.$$

for all $\alpha,\beta \in \Delta$.

There is a natural isomorphism $\psi:WL\to$Pic($G/B$).
The map is defined as follows:  for each $\omega \in WL^{+}$
there exists a representation $\rho: G \to GL(V)$ with highest weight 
$\omega$. Let the highest weight vector be $v\in V$. Then
there is a map $G/B \to Orb(\Bbb{C}v)\subset \Bbb{P}(V)$ where 
$Orb(\Bbb{C}v$) is the orbit of the line $\Bbb{C}v$. The map then takes 
$\omega$ to pull back of $\mathcal{O}(1)$ by the map above.

The subset $WL_{P}$ of weights that descend to line bundles on $G/P$ are just those elements
$\omega$ which satisfy the property

if $\frg_{\alpha}\bigoplus\frg_{-\alpha}\subset\frp$, then $\omega(H_{\alpha})=0.$

The second Homology group $H_2(G/P,\Bbb{Z})$ can be naturally considered as 
the dual $Hom(WL_{P},\Bbb{Z})$ by Poincare duality. Note that the homology
class of $f_*([C])$ where $f:\Bbb{P}^1\to G/P$ corresponds to the map
$WL_P \to \Bbb{Z}$ obtained by taking 'degree of the pullback bundle'.

Next we describe the first Chern classes of the tangent bundles of $G/B$ and
$G/P$. 
\begin{enumerate}
\item First Chern class of $T_{G/B}$: this is $\psi(\sum_{\alpha \in R^{+}}
\alpha).$
\item First Chern class of $T_{G/P}$: this is $\psi(\sum_{\alpha \in R,\frg_{\alpha}\subset \frp}\alpha).$
 \end{enumerate}

\subsection{Cell decomposition and cohomology of $G/P$}

Let $P\subset B$ be a parabolic subgroup. Let $\Delta_P$ be the set of roots $\alpha$ such that $\frg_{\alpha }\bigoplus \frg_{-\alpha}\subset \frp.$ Let $W_{P}$ be the subgroup of the Weyl group generated by the reflections corresponding to elements of $\Delta$.

\begin{te}(Bruhat decomposition:) $G/P$ is a disjoint union of the sets
$\Lambda_{w}$ for $w \in W/W_{P}$ where $\Lambda_w$ is defined to be $BwP
\subset G/P$. Let $X_w$ be the closure of $\Lambda_w$. The codimension of $X_w$ is the cardinality of the set

$$\left|\{ \alpha \in R:\frg_{\alpha}\not\subset \frp, w(\alpha)\not\in R^{+}\}\right|.$$

\end{te}
\begin{proof}Standard, see \cite{bor}.
\end{proof}
It is known that the subvarieties $X_w$ generate the cohomology (additively)
of 
$G/P$.Finally recall the definition of 'relative position' $[g_1,g_2]$ of two elements $g_1,g_2 \in G$. We define $w=[g_1,g_2]$ to be the unique element of the Weyl group so that there exist $b_1,b_2 \in B$ so that
$g_1 =g_2 b_1 w b_2$ (Bruhat decomposition).

Note the following three properties:
\begin{enumerate}
\item $[g_1,g_2] = [gg_1,gg_2]$ for $g_1, g_2 \in G.$
\item $[g_1 b,g_2]=[g_1,g_2 b]= [g_1,g_2]$ for $b\in B.$
\item $h \in g \Lambda_w$ if and only if $[h,g]=w$.
\end{enumerate}

Analogous definition of relative position can be made of $[g_1,g_2]$ where
this takes values in $W/W_{P},g_1\in G/P$ and $g_2\in G$.

We need one final lemma which relates the codimensions of $X_w$ and $X_{w_c w}$
where $w_c$ is an element of the Weyl group constructed out of a central element $c$ as in the previous section

\begin{nle} Let $c\in \text{center}(G), w \in W$ with representative $x_{\alpha}\in \overline{C}$. Then $\text{codim}(X_{w_c w})-\text{codim}(X_{w})$
is equal to

$$\sum_{\beta \not\in \Delta} w\beta(x_{\alpha}).$$
\end{nle}

\begin{proof}

The quantity we are interested in is
$$\left|\{ \alpha \not\in \Delta_P: w(\alpha)\in R^{+}\}\right|- \left|\{ \alpha \not\in \Delta_P:(w_c w)(\alpha)\in R^{+}\}\right|.$$

To evaluate the second quantity(using lemma 1 of 2.4), divide into two cases namely

$$\left|\{ \alpha \not\in \Delta_{P}: w(\alpha)\in R^{+},w\beta(x_{\alpha})=0\}\right|$$

and

$$\left|\{ \alpha \not\in \Delta_{P}: w(\alpha)\in R^{-},w\beta(x_{\alpha})=-1\}\right|.$$

Therefore the quantity we are interested in becomes

$$\left|\{ \alpha \not\in \Delta_{P}: w(\alpha)\in R^{+},w\beta(x_{\alpha})=1\}\right|- \left|\{ \alpha \not\in \Delta_{P}: w(\alpha)\in R^{-},w\beta(x_{\alpha})=-1\}\right|.$$

and that is what is displayed in the statement of the lemma. Note that our
hypotheses imply that if $\beta\in R$, then $\beta(x_{\alpha})$ is in the set
$\{-1,0,1\}$.

\end{proof}

\subsection{Space of maps and Gromov-Witten invariants}

Let $X \in H_{2}(G/B; \Bbb{Z})$.

Let $M_{X} = $ space of maps $f:\Bbb{P}^1 \to G/B$ so that
$f_{*}([\Bbb{P}^{1}]) =X$. It is known that $M_{X}$ can be given the structure of
a smooth quasi-projective variety of dimension $d_X$ where

$d_X = c_1(T_{G/B})\cap X$ + dim($G/B).$

Of central interest to us in this paper are Gromov-Witten invariants. Recall
that we have fixed three points $p_1,p_2,p_3$ on $\Bbb{P}^1$ (which we usually take to be $0,\infty,1$).

\begin{defi} Let $Z \in H_{2}(G/B; \Bbb{Z})$; $w_1,w_2,w_3 \in W$
. Then, 

$<X_{w_1},X_{w_2},X_{w_3}>_{Z}$ is defined to be the number of
maps (0 if infinite) $f \in M_Z$ so that $f(p_i) \in g_i X_{w_i}; i=1,2,3$
where $g_i$ are 'general' points of $G$.
\end{defi} 
Note that the invariant above is 0 unless
$$\sum \text{Codim}(X_{w_i})= c_1(T_{G/B})\cap Z + \text{dim}(G/B).$$

\section{The transformation formula}

\begin{te} 

Let $Z \in H_{2}(G/P; \Bbb{Z})$; $u_1,u_2,u_3 \in W$
,$c_1,c_2\in \text{Center}(G)$. Then,

$$<X_{u_1},X_{u_2},X_{u_3}>_{Z}=<X_{w_{c_1}u_1},X_{w_{c_2}u_2},X_{u_3}>_{Z'}.$$

where $Z'$ as an element of $Hom(WL_P,\Bbb{Z})$ is given by

$$z'(\gamma)= z(\gamma)-\gamma(u_1^{-1} x_1-x_1) -\gamma(u_2^{-1}x_2 - w_2 x_2).$$
\end{te}

\begin{proof} Let us first check that the codimension condition
$$\sum \text{codim}(X_{u_i})= c_1(T_{G/P})\cap Z + \text{dim}(G/P).$$

for the left hand side is the same as that for the RHS.

Recall that if $c\in \text{center}(G), w \in W$ with representative $x_{\alpha}\in \overline{C}$. Then $\text{codim}(X_{w_c w})-\text{codim}(X_{w})$
is equal to

$$\sum_{\beta:\frg_{\beta}\not\subset p} w\beta(x_{\alpha}).$$

so we have to verify

$$\sum_{\beta:\frg_{\beta}\not\subset p} (u_1\beta(x_1)+u_2\beta(x_2))$$

equals

$$ \gamma(Z')-\gamma(Z)$$ where 
$$\gamma =-\sum_{\beta:\frg_{\beta}\not\subset p}\beta.$$ 

which has been proved before [3.2, lemma 4].

To start working towards the transformation formulas, fix $u_1,u_2,u_3 \in W$,
$g_1,g_2,g_3\in G$ 'elements in general position'. And $c_1,c_2$ central elements in $G$ with $c_1 c_2 =1$ (as in section 2.5). Recall the notations of
section 2.5.

Pick an element $k\in (P_1 g_1^{-1})\cap (w_{c_2}P_2 g_2^{-1})$. There exists such a $k$ because of corollary 1 in section 2.5. Consider the map

$\phi: G/P \to G/P$ given by left multiplication by $k$. Let  $f \in M_Z$,
we claim

$f(p_i) \in g_i X_{w_i}; i=1,2,3$ if and only if

Setting $g=\phi f ,g(p_i) \in kg_i X_{w_i}; i=1,2,3$, this claim is obvious
but we also have $k(g_1) \in P_1$ and $kg_2\in Ad(w_{c_2})P_2$
. So we might as well assume $g_1 \in  P_1$ and $g_2 \in w_{c_2}P_2$.

Now suppose $f: \Bbb{P}^1 \to G/P$ then , let $s$ be the map $\Bbb{P}^1 \to
G/P$ given by:
(to simplify notation: let 
\begin{enumerate}
\item $w_1$ and $w_2$ denote $w_{c_1}$ and $w_{c_2}$ respectively.
\item $x_1$ and $x_2$ denote $x_{\alpha}$ and $x_{\beta}$ respectively.
\item for $t\in T$, $z \in \Bbb{C}$, let
$z^{t}$ denote $exp(ln(z)t)$ (with the indeterminacies.)
\end{enumerate}

$$s(z)= z^{ x_1} f(z).$$

Note that the indeterminacy of $s$ is always {\em central} so as a map to
$G/P$ it is well defined (on the complement  of $\{0,1,\infty\}$). And we
can extend this to the whole of $\Bbb{P}^1$, because all the functions involved are of bounded growth.

We have to study the effect on the degrees and also on the 'positions' of
$s(1),s(1),s(\infty)$.

{\em Position of s(0):}\\

We have assumed that $g_1$ is in $P_1$. And $P_1 =Q_1B$. So let $g_1$
= $q_1b$. We claim that the element

$g_1'= q_1 w_1^{-1}$ is well defined in $G/P$ (independent of choices). 
we need that if $q_1'=q_1 b$, then$q_1 w_1^{-1}$ and $q_1 b w_1^{-1}$ give the same point
in $G/P$. That is, $w_1 b w_1^{-1} \in B$ if $ b \in Q_1 \cap B$. But this is clear from lemma 1 of section 2.4.

We claim:
$[s(0),g_1']=w_1[f(0),g_1]$.
That is
$[s(0),q_1 w_1^{-1}] =w_1^{-1}[f(0),q_1]$.

Let $f(0)=q_1 b_1 w b_2$ and $f=n(z)f(0)$ where $n(0)=1$. we therefore need to compare

$[lim_{z\rightarrow 0}z^{x_1}n(z)f(0),g_1']$ with

$[f(0),g_1]$.

Or,if we set $h=q_1^{-1}f(0) $,

we want to relate

$[lim_{z\rightarrow 0}(Ad(q_1^{-1})z^{ x_1}) (Ad(q_1^{-1})n(z))h,w^{-1}]$

with $[h,1]$

It is easy to see that $Ad(q_1^{-1})z^{ x_1} = z^{ x_1}$.
Setting $r(z) = Ad(q_1^{-1})n(z)$, we want then to relate

$[lim_{z\rightarrow 0}z^{ x_1}r(z)h,w^{-1}]$ with $[h,1]$ where $r(0)=1$. For this we need

\begin{nle}
If $d(z)$ is a holomorphic map to $G$, with $d(0) \in B$ then
$k =lim_{z\rightarrow 0}Ad(z^{x_1})d(z)$ exists with $w_1 k w_1^{-1} \in B$.
\end{nle}
\begin{proof} $G$ is generated by the one parameter groups $G_{\alpha}$ for
$\alpha \in R$ and T. These groups are isomorphic to $\Bbb{C}$, and with an action (Ad) of the torus with
$Ad(t)u = \alpha(t)u$.

If $d(z) \in  G_{\alpha}=\Bbb{C}$ given by $d(z)=z^m$ then
 
$Ad(z^{x_1})d(z) = z^{\alpha(x_1)}d(z)$. Hence to verify the lemma
we need

\begin{enumerate}
\item If $\alpha \in R^{+}$ then $\alpha(x_1)= 0$ implies $w_1(\alpha)$ is a positive root, which is known. If $\alpha(x_1)= 1$ then $k=1$.
\item If $\alpha \in R^{-}$ with $\alpha(x_1)= 0$ then clearly $k=1$.
\item If $\alpha \in R^{-}$ with $\alpha(x_1)= -1$ then clearly $k$ exists
and $w_1(\alpha)$ is positive.
\end{enumerate}
\end{proof}

{\em Position of $s(\infty), s(1)$:}
Note that we have chosen $p_2=\infty$ in order to simplify the notation.

Note that $g_2 \in w_2 P_2$. Consider the map $\psi: G/P \to G/P$ by left
multiplication by $w_2^{-1}.$ Write $g_2 = w_2 q_2 b$ for $q_2 \in Q_2$
 and $b \in B$. Let $g_2' = w_2 q_2 w_2^{-1}b$ as before the 
$G/P$ class of $g_2'$ is well-defined. It is then easy to see that
$w_2[f(\infty),g_2]=[s(\infty),g_2']$. For this is is enough to notice that
$z^{x_1}= (\frac{1}{z})^{w_2 x_2}$.

Finally let $g_3' =g_3$ it is clear that $[s(1),g_3']=[f(1),g_3]$.

Next we have to compute the homology class $s_*([\Bbb{P}^1])$. 

\begin{te}
If $z$ is the element in $Hom(WL_P,\Bbb{P}^1)$ corresponding to
$f$ then the element $z'$ corresponding to $s$ is the element

$$z'(\gamma)= z(\gamma)-\gamma(u_1^{-1} x_1-x_1) -\gamma(u_2^{-1}x_2 - w_2 x_2)$$
$$ = z(\gamma)-u_1 \gamma(x_1) -u_2 \gamma(x_2).$$
\end{te}
\begin{proof}
It is enough to prove this in the case $\gamma$ positive and integral.
Let $L$ be the line bundle on $G/P$ corresponding to $\gamma$. We construct the line bundle corresponding to $\gamma$ in a different (equivalent) manner first.
First extend $\gamma$ to a map $\Gamma:P\to \Bbb{C}^*$. then construct
the total space of $L$ as $G\times\Bbb{C} /R $ where $R$ is the equivalence relation $(g,v)=(gp,\Gamma(p)v)$ for $p\in P$. 

$f,s$ give two line bundles $L_f = f^{*}(L),L_s=s^{*}(L).$ At a point other
than 0,$\infty$, construct the map $\psi:L_f \to L_s$, by
$(\breve{f}(z),1) $ to $(z^{x_1}\breve{f}(z),\Gamma(z^{x_1}))$ where
$\breve{f}(z)$ is a local lifting of $f$ to a map $ \to G$, and where the
same determination of $z^{x_1}$ in both $z^{x_1}\breve{f}(z)$ and in $\Gamma(z^{x_1})$.

It is immediate to see that $\psi$ is an isomorphism of bundles outside
of $\{0,\infty\}$. Let us analyze this map first at $z=0$. Lift $f$ to a map
$\breve{f}$ to $G$. Then $(\breve{f},1)$ is a local section of $f$ and
this is mapped by $\psi$ to $(z^{x_1}\breve{f},\Gamma(z^{x_1}))$ a meromorphic
section of $L_s$. To complete the analysis we have to display a generating
section of $L_s$. Let $\breve{f}(z)=q_1b(z)u_1 p$ where (recall $g_1=q_1 b_1 p \in P$) $b(0)\in B$. Now $z^{x_1}q_1=q_1 z^{x_1}$ and $z^{x_1}b(z)z^{-x_1}$  is holomorphic at $z=0$. We therefore find that $(z^{x_1}q_1 b(z)z^{-x_1}u_1 p_1,1)$ is a holomorphic section of $L_s$. Therefore the contribution at
$z=0$ to $deg(L_s)-deg(L_f)$ is $\gamma(x_1-u_1^{-1}x_1)$.

The calculation at $\infty$ is similar and we arrive at the equation 
in the statement.

\end{proof}

Now consider the map which takes a map $f:\Bbb{P}^1\to G/P$ to the map
$s$ as above. We have seen that if we choose generic $g_i$
to compute the left hand side, then the $s's$ correspond to the right hand side computed
with respect to $g_i'$. The $g_i'$ depend only on the $g_i$ and the central elements chosen. So computed w.r.t $g_i'$ the right handside is a finite
number and the codimension computation therefore gives us an inequality.

$$<X_{u_1},X_{u_2},X_{u_3}>_{Z}\leq<X_{w_{c_1}u_1},X_{w_{c_2}u_2},X_{u_3}>_{Z'}.$$

Now apply the reasoning again, this time with $c_2,c_1$, to get the other
inequality.

\end{proof}

\begin{cl} 

Let $Z \in H_{2}(G/P; \Bbb{Z})$; $u_1,\dots, u_s \in W$
,$c_1,\dots,c_s \in \text{Center}(G)$,\\$ c_1 c_2 \dots c_s =1$. Then,

$$<X_{u_1},\dots,X_{u_3}>_{Z}=<X_{w_{c_1}u_1},\dots,X_{w_{c_s} u_s}>_{Z'}.$$

where $Z'$ as an element of $Hom(WL_P,\Bbb{Z})$ is given by

$$z'(\gamma)= z(\gamma)-\gamma(u_1^{-1} x_1) -\gamma(u_2^{-1}x_2)-\dots - \gamma(u_s^{-1} x_s).$$

($x_{k}\in \overline{C}$ is the representative of $c_k$.)

\end{cl}
\begin{proof} 
Let us do the case $s=3$, the general case is similiar. We write down the
transformation formulas (as in the theorem) for $c_1,{c_1}^{-1}$, and then transform this on the second and third 'variables' by $c_1c_2, (c_1c_2)^{-1}$, it is clear that $(c_1c_2)^{-1} =c_3$. We just have to verify that the formula for
$Z'$ is the one above. We leave this to the reader.
\end{proof}

\section{Reformulation}
Let $P$ be a standard parabolic, $\Sigma$ the set $\{\alpha \in \Delta:\frg_{-\alpha} \not\subset \frp\}$. It is clear that $\omega_{\sigma}:\sigma \in \Sigma$ gives a basis for $WL_P$. Introduce variables
$q_{\sigma}$ for $\sigma \in \Sigma$.

We need the following simple fact before we can describe the quantum cohomology of $G/P$.

{\em Duals: 

}The Poincare dual of the class $X_{w}$ is the class $X_{w_{0} w}$ where $w_0$
is the unique element in the Weyl group so that

$$B\cap w_0B{w_0}^{-1} = T. $$

\begin{defi} Define 

\begin{enumerate}
\item $X_{u_1}\star X_{u_2}$\\$ = \sum_{u\in W/W_P,Z \in Hom(WL_P,\Bbb{Z})}(\prod_{\sigma\in \Sigma}{q_{\sigma}}^{Z(\omega_{\sigma})})<X_{u_1},X_{u_2},X_{u}>_Z X_{w_0 u}.$

\item $QH(G/P)= H^*(G/P,\Bbb{C}) \bigotimes \Bbb{C}[q_\sigma: \sigma \in \Sigma]$ with the product given above.

\item For $c$ in center of $G$ with $w_c$ the associated Weyl group element,
let $T_c:QH(G/P)\to QH(G/P)$ by $T_c(X_w)=(\prod_{\sigma}{q_{\sigma}}^{\omega_\sigma(x_1-w^{-1}x_1 )}X_{w_c w}$.
\end{enumerate}

\end{defi}

Let us now try to compare $T_c(X_{u_1} \star X_{u_2})$ to $T(X_{u_1})\star X_{u_2}$

\begin{nle}
\begin{enumerate}
\item $T_{c_1}T_{c_2}= (\prod_{\sigma}{q_{\sigma}}^{\omega_\sigma(w_1^{-1} x_2 - x_2)})T_{c_1 c_2} 
$ as operators.

\item $T_c (x \star y) = T_c (x) \star y$. 

\item $T_1$ = multiplication by   $1.$

\end{enumerate}
\end{nle}
\section{The $SL_N$ case}

Let us look at $Gr(r,n)$. Here we have simple roots $L_{i}-L_{i+1}$ for
$i=n-1,\dots,1$. For $Gr(r,n)$, $\Sigma$ from the previous section is
$L_{r}-L_{r+1}$. The center is cyclic group of order $n$ generated by
the diagonal matrix $\Theta$ with entries $\zeta$ where $\zeta =e^{\frac{2\pi i }{n}}$.

The element in $\bar{C}$ corresponding to $\Theta^k, k=1,\dots n-1$ is
$(\frac{k}{n},\dots,\frac{k}{n},\frac{k}{n}-1,\frac{k}{n}-1)$ where there
are $(n-k)$ $ \frac{k}{n}$'s. The element of the Weyl group corresponding to
$\Theta^k$ is just 'subtract k modulo n,replacing zeroes by n' in the standard
representation of the Weyl group as a permutation group. Using these we
can give a more explicit form of the transformation formulas.

\begin{defi}
Let $I = \{ i_{1} < i_{2} < \dots < i_{r} \}$. Let $F^{\sssize{\bullet}}$ be a complete flag in a n-dimensional vector space E. Now let
$\Omega_{I}(F^{\sssize{\bullet}}) = \{ L \in Gr(r,E) \mid \text{dim}(L \cap
F^{i_{t}}) \geq t$ for $1\leq t \leq r \}$. We denote the cohomology class of
this subvariety by $\sigma(I)$. The codimension of this subvariety is the
number of pairs $(j,i)$ with $j \notin I$ and $i\in I$ and $j>i$.
\end{defi}

        The Gromov-Witten invariants in the Grassmannian case also have an interpretation in terms of
vector bundles on $\Bbb{P}^{1}$.  Let $V = \mathcal{O}^{n}$ be a vector bundle on 
$\Bbb{P}^{1}$ . We have a universal sequence of vector bundles on $Gr(r,n)$

$$
0 \to \mathcal{S} \to \mathcal{O}^{n} \to \mathcal{Q} \to 0
$$ 

 where $ \mathcal{S}$ is the universal subbundle of rank $r$ and $\mathcal{Q}$ the quotient. It is now easy to verify that degree $d$ maps  $\rho  :\Bbb{P}^{1} \to Gr(r,n)$ are in $1$-$1$ correspondence with subbundles of rank $r$ and degree $-d$ of V by pulling back the universal sequence via the map $\rho$. Also, the image of point $p_{i}$ under this map is exactly the fiber of this subbundle at $p_{i}$. It is useful to fix an $n$-dimensional space $T$ and  identify all fibers of the bundle V with $T$. To obtain the other direction of this correspondence note that subbundles $\mathcal{S}$ correspond to a family of $r$ dimensional subspaces of T(over $\Bbb{P}^1$).

        The number defined above, therefore counts the number (zero if infinite) of subbundles $E$ of V of degree $-d$ and rank $r$ such that the fiber $E_{p_{i}}$ as a subset of T lies in the Schubert variety $\Omega_{I}(F^{\sssize{\bullet}}_{p_i})$.

\subsection{The Transformation property}

\begin{te}
Let $I_1,\ldots, I_s$ be subsets of $\{1,..n\}$ of cardinality $r$ each. Let $n_1,\ldots,n_s$
be natural numbers summing to $n$.
Define $J_i = I_i-n_i$ mod $n$. That is subtract $n_i$ from the numbers in $I_i$
and reduce them mod $n$ and replace all $0$'s by $n$.
       define $d_i$ = number of elements in $I_i$ which are less than or equal to $n_i$.

then
$$<\sigma(I_1),\ldots,\sigma(I_s)>_d = <\sigma(J_1),\ldots,\sigma(J_s)>_{d+r-\sum d_i}.$$
\end{te}

Remark: This property has been noted and proved in
\cite{wood} as a consequence of the Schubert Calculus established in \cite{bert}. The proof here is geometric and independent of \cite{bert}. It is essentially the same proof as that of the transformation formulas of the previous section.

\begin{proof}
First we verify that the codimension conditions on both sides are the same.
That is:
$$\sum_{i=1}^{s} \text{codim}(\sigma(I_i)) = nd +r(n-r). $$

is same as the condition
$$\sum_{i=1}^{s} \text{codim}(\sigma(J_i)) = n(d+r -\sum{d_i}) +r(n-r). $$ 

This follows easily from the observation:

$$\text{codim}(\sigma(J_i)) = \text{codim} (\sigma(I_i)) + (n_i- d_i)r - (n-r)d_i.$$

Now fix $s$ general flags on $V=C^n$: $F^{\sssize{\bullet}}_{p_{i}}, i=1,..s$, as well as $s$ points $p_1,\dots,p_s$ on $\Bbb{P}^1$.
Then  $<\sigma(I_1),\ldots,\sigma(I_s)>_d$ is the number of subbundles
        (zero if infinite) of subbundles $\mathcal{E}$ of $\mathcal{V}=V\bigotimes_{\Bbb{C}}\mathcal{O}$ of degree $-d$ and rank $r$ such that the fiber $\mathcal{E}_{p_{i}}$ as a subset of $C^n$ lies in the Schubert variety $\Omega_{I_i}(F^{\sssize{\bullet}}_{p_i})$. Let $V_i = F^{n_i}_{p_i}$. We have
from genericity of the flags, an equality

$$\bigoplus V_i \mapsto C^n.$$

We define a new bundle on $\Bbb{P}^1$ as follows ($q$ is a new point on $\Bbb{P}^1$).
\begin{defi}
$\mathcal{V}' = \bigoplus (V_i\bigotimes\mathcal{O}(p_i-q))$.
\end{defi}

Note that we are given an isomorphism $\mathcal{V}\to\mathcal{V}'$ over
the open set $U = \Bbb{P}^1 - \{p_1,\dots,p_s,q\}$ and also that
$\mathcal{V}'$ is isomorphic to $\mathcal{O}^{n}$. From the theorem
below we know that we have a one-one correspondence between subbundles
$\mathcal{E}$ of $\mathcal{V}$ and subbundles $\mathcal{E}'$ of $\mathcal{V}'$.
We also have induced flags  $F'^{\sssize{\bullet}}_{p_i}$
on the fibers of $\mathcal{V}'$ at the points $p_i$,
so that if fiber $\mathcal{E}_{p_{i}}$ is in the   Schubert variety $\Omega_{I_i}(F^{\sssize{\bullet}}_{p_i})$, then the  fiber $\mathcal{E'}_{p_{i}}$ is in the Schubert variety $\Omega_{J_i}(F'^{\sssize{\bullet}}_{p_i})$. The theorem below also tells us that in this correspondence degree of $\mathcal{E}'$ is
equal to degree of $\mathcal{E} -r + \sum d_i$.

This finishes the proof but we have to deal with genericity questions. A more
refined approach can directly show that the induced flags on $\mathcal{V}'$
are generic too. But we wish to avoid this line of argument here. Instead
we note that this argument proves that (with the codimension computation) that

$$<\sigma(I_1),\ldots,\sigma(I_s)>_d        \leq
 <\sigma(J_1),\ldots,\sigma(J_s)>_{d+r-\sum d_i}.$$

We could then reverse this construction (or perform in many times) to
get the other inequality. This proves that there are no intersections at `
infinity' and also transversality without invoking the theory of Quot schemes.
\end{proof}

\begin{te}(Local theory) Let $C$ be a smooth curve $p$ a point on it, $t$ a uniformising parameter at $p$. Let $\mathcal{V}$ be a vector bundle on $C$,
$V=$ fiber  $\mathcal{V}_{p}$. Also suppose that we are given a complete flag on $V$: $V_1\subset V_2 \subset \dots \subset V_n = V$. Define
$\mathcal{V}'_k = \{$ meromorphic sections $s$ of $\mathcal{V}$ which are holomorphic sections of $\mathcal{V}$ outside of $p$ and such that $ts$ extends to give a section of $\mathcal{V}$ near $p$ with fiber at $p$ in $V_k \}$. 
\begin{enumerate} 

\item $\mathcal{V} \subset \mathcal{V}'_k$ with quotient supported at $p$ of dimension $k$.
\item There is a one-one correspondence between subbundles $\mathcal{E}$ of $
\mathcal{V}$ and subbundles $\mathcal{E}'_k$ of $\mathcal{V}'_k$. In this correspondence the quotient $\frac{\mathcal{E}'_k}{\mathcal{E}}$ is supported at
$p$ and has dimension $\text{dim}(\mathcal{E}_p \cap V_k)$.
\item We have a sequence of inclusions

$$
    t(\mathcal{V}'_k) \subset t(\mathcal{V}'_{k+1})\dots\mathcal{V}\subset \mathcal{V}'_1\dots \subset \mathcal{V}'_k.$$

which gives a complete flag on the fiber $(\mathcal{V}'_k)_p$.

\item In the correspondence on the subbundles if fiber $\mathcal{E}_{p_{i}}$ is in the   Schubert variety $\Omega_{I}(F^{\sssize{\bullet}}_{p_i})$, then the  fiber $\mathcal{E_k'}_{p}$ is in the Schubert variety $\Omega_{I-k}(F'^{\sssize{\bullet}}_{p})$. 

\end{enumerate}
\end{te}
\begin{proof} All fairly obvious and proofs can be found in the appendix to
\cite{be}.
\end{proof}

\subsection{Quantum Schubert calculus}

The objective of this section is to show how Pieri's formula is a easy consequence of the relations of the previous section. Namely for the intersections
in the Pieri formula, the relations reduce to $d=0$ case which are classically known by induction.

As in the classical approach of Hodge-Pedoe see \cite{gh}, one can then
prove Giambelli formula from Pieri. In \cite{bert} Pieri is deduced from
Giambelli.

\begin{defi} Define a map $T:QH(Gr(r,n))\mapsto QH(Gr(r,n))$ by the rule
$T(\sigma(I))= q^{d_1}(\sigma(I-1))$ where $d_1$ = number of elements in $I$ which are less than or equal to 1. Notice that if $k<=n$, $T^k(\sigma(I))= q^{d_k}\sigma(I-k)$ where $d_k$ = number of elements in $I$ which are less than or equal to $k$. This $T$ is essentially $T_{\Theta}$ of the previous section.
\end{defi}

\begin{te}(reformulation of transformation property)
The transformation property is equivalent to the property
$$T(\sigma(I)\star\sigma(J))=T(\sigma(I))\star \sigma(J).$$
\end{te}

The proof of Pieri to follow is pure algebra beyond this point. We first note
the following
\begin{nle}
Let $\sigma(I)$ have codimension less than or equal to $n-1$, then if
$\sigma(K)$ appears in $\sigma(I)\star\sigma(J)$ with a q coefficient $>=1$, then  then there exists $k$ so that $\sigma(K-k)$ appears in  $\sigma(I)\star\sigma(J-k)$ with q degree $=0$.
\end{nle}
\begin{proof}
 Suppose not, choose $k$ so that the q degree is minimized and equals $d$. Let $J'=J-k$, and
$K'=K-k$. since we cannot minimize the q degree further, application
of the operation $T^l$ tells us that if $d_l$ = number of elements in
$J'$ which are less than or equal to $l$, and $c_l$ = similar number for $K'$
then $d_l \leq c_l$. This clearly implies that $\text{codim}(J')\leq \text{codim}(K')$.
But we also have 

$$\text{codim}(I) + \text{codim}(J') = nd + \text{codim}(K')$$ this yields a contradiction immediately if $d>=1$.
\end{proof}

Pieri formula is usually written in cohomological notation: for this we make the following definitions.
\begin{defi} If $I=\{i_1<\dots<i_r\}$ is a subset of $\{1,\dots n\}$ then
define $a(I,k)= n-r+k-i_k$ for $k=1,\dots, r$.
\end{defi}
\begin{defi}(special schubert cells)
If $a<=n-r$ define $\sigma_a = \sigma(I_a)$ where 

$$I_a =\{n-r+1-a, n-r+2,\dots n-r\}.$$
\end{defi}

\begin{te}
$$\sigma_a\star \sigma(I) =\sum_K \sigma(K) + q\sum_L \sigma(L).$$ 
where the $K$ sum is over all $K$ satisfying

$$n-r\geq a(K,1)\geq a(I,1)\geq a(K,2)\geq \dots\geq a(K,r)\geq a(I,r)$$

and $$\text{codim}(I) + a = \text{codim}(K).$$

and the $L$ sum is over all $L$ satisfying

$$ a(I,1)-1\geq a(L,1)\geq a(I,2)-1\geq \dots\geq a(I,r)-1\geq a(L,r)\geq 0$$

and $$\text{codim}(I) + a = \text{codim}(L) + n .$$

Notice that there are no $L$ terms if $a(I,r)=0$.

\end{te}

\begin{proof}
The statement about the $K$ terms is classical \cite{gh}. Let us first show that there are no terms with $q^2$ and higher. From the previous lemma there
is then $L'$ and  and $I'$ so that $I=I'-k$ and $L=L'-k$, so that
\begin{itemize}
\item $\sigma(L')$ appears with $q$ degree =0 in $\sigma_a\star\sigma(I')$
\item $d_k <=c_k -2$ where $d_k$ = number of elements in $I' \leq k$, and
$c_k$ = number of elements in $L' \leq k$.
\end{itemize}

Let $d_k = j$. Then $a(L',j+2)\leq a(I',j+1)$ tells us that if
$L'=\{l'_1\leq\dots\leq l'_r\}$ and $I'=\{i'_1\leq\dots\leq i'_r\}$
then

$$n-r+(j+1)-l'_{j+2} \leq n-r+j-i'_{j+1}$$ 

or what is the same

$i'_{j+1} +1 \leq l'_{j+2}$ hence $l'_{j+2} \geq i'_{j+1} +1$. But
$i'_{j+1}>k$ , therefore $l'_{j+2} > k+1$ which is in direct contradiction to
$d_k <=c_k -2.$

We now deal with the $q^1$ terms. We find $L'$ and $I'$ as before
so that $d_k = c_k -1$. Let $d_k =j$. $c_k =j+1$. Then contemplation
of $I=I'-k$ and $J= J'-k$ gives the result. For example let us only verify
$a(I,r)-1\geq a(L,r)$. It is clear that
$l_r = l'_{j+1} -k +n$, and $i_r = i'_j -k +n$.

so we need to verify
$$n-r +r -l_r \leq n-r +r -i_r - 1$$
or that
$$i_r+1 \leq l_r$$

or that
$$ i'_j +1 \leq l'_{j+1}.$$ which is just
$$a(I',j)\geq a(L',j+1)$$ which is already known classically. 

\end{proof}

The tediousness of the above proof is made up by the simplicity of the proof
of Giambelli formula. The essential part of which is the the following formula.
Because of our Homological -cohomological notation problems we need to make one definition.
\begin{defi}
Let $\breve{a} =a_1\geq\dots\geq a_r$ be given, then define
$I(\breve{a})$ by $i_l = n-r+l -a_l$. If this is not a subset of $\{1,..n\}$
define $\sigma(I(\breve{a}))$ to be zero. Also denote  $\sigma(I(\breve{a}))$
by $\sigma_{a_1,\dots,a_r}$ and ignore 0's in the subscript if they appear.
\end{defi}
\begin{te}

$$(-1)^d \sigma_{a_1,\dots a_d} =\sum_{j=1}^{d} (-1)^j \sigma_{a_1,\dots,a_{j-1},a_{j+1}-1,\dots,a_d-1}\star \sigma_{a_j+d-j}.$$
\end{te}

\begin{proof}
Note that on the lhs the length of the string $\breve{a'}= \{a_1,\dots,a_{j-1},a_{j+1}-1,\dots,a_d-1 \}$ is less than r. So
$a'_r =0$. Hence no $q^1$ terms are produced via Pieri, and the formula in $q$ degree 0 is known  classically \cite{gh} page 205. Hence there is really nothing
(new) to prove. Iteration of this gives the Giambelli formula as in \cite{gh}.
\end{proof}

Fulton and Woodward have proved a theorem on the smallest power of q in the quantum product of Schubert subvarieties in the case of $G/P$, P maximal parabolic. We prove the Grassmann case of their theorem in a slightly strengthened form.
\begin{te}(Fulton-Woodward) The smallest power of q appearing in
$\sigma(I)\star\sigma(J)$ is the number
$d=max\{d_i+d'_j-r:i+j=n\}$ where $d_i$ = number of elements in $I$ which are
less than or equal to $i$, $d'_j$ = number of elements in $J$ which are
less than or equal to $j$. Moreover, if the max is achieved for $i,j:i+j=n$, then

$$\sigma(I)\star\sigma(J)= q^d(\sigma(I-i)\cup\sigma(J-j)) + \text{higher order terms}.$$

where the cup product $\cup$ on the right hand side is the cup product in usual cohomology.

\end{te}

{\em Remark:} The `strengthened' part refers to identification of the lowest
order terms which is curiously a product in the ordinary cohomology. This may not be true for all $G/P$'s. Also the associativity of quantum cohomology has not been used so far!

\begin{proof} We know $T^n =$ multiplication by $q^r$, therefore

$$q^r(\sigma(I)\star \sigma(J))= T^n(\sigma(I)\star \sigma(J)).$$

$$= T^i(\sigma(I)\star T^j(\sigma(J)$$ when $i+j =n$,

$$=q^{d_i+d'_j} \sigma(I-i)\star \sigma(J-j).$$

Therefore, if $i+j=n$,
$$\sigma(I)\star \sigma(J) = q^{d_i+d'_j-r} \sigma(I-i)\star \sigma(J-j).$$

Suppose that we choose $i,j$(with sum $=n$)  so as to maximize $d_i+d'_j-r$
then clearly for $I'=I-i$  and $j'=J-i$, we have $c_k +c'_{n-k}\leq r$, for
any $k$ (here $c_k$= number of elements in $I'$ which are
less than or equal to $k$, $c'_{n-k}$ = number of elements in $J'$ which are
less than or equal to $n-k$. ) this tells us that the dual of $\sigma(I')$
is contained in $\sigma(J')$, for a choice of flags. Hence by
 Kleiman's Bertini theorem, we have $\sigma(I')\cap\sigma(J')\neq 0$.

\end{proof}

\end{document}